\documentclass[12pt]{article}       

\usepackage{amsmath,amsthm}
\usepackage{amssymb}
\usepackage{latexsym}
\begin{document}

\newcommand{\comment}[1]{}    
\newcommand{\hs}{\enspace}
\newcommand{\hhs}{\thinspace}
\newcommand{\real}{\ifmmode {\rm R} \else ${\rm R}$ \fi}

\newtheorem{theorem}{Theorem}
\newtheorem{lemma}[theorem]{Lemma}         
\newtheorem{corollary}[theorem]{Corollary}
\newtheorem{definition}[theorem]{Definition}
\newtheorem{claim}[theorem]{Claim}
\newtheorem{conjecture}[theorem]{Conjecture}
\newtheorem{proposition}[theorem]{Proposition}

\newtheorem{remark}[theorem]{Remark}



\def\edge{\leftrightarrow}
\def\noedge{\not\leftrightarrow}
\def\twoedge{\Leftrightarrow}
\def\to{\rightarrow}
\def\Hrl{H^{(r)}_{l+1}}
\def\Krl{{\cal K}^{(r)}_l}
\def\Krl{{\cal K}^{(r)}_{l+1}}
\def\cF{{\cal F}}
\def\cG{{\cal G}}
\def\cH{{\cal H}}
\def\e{\varepsilon}
\def\bF{
{\cal {\bf F}}}
\def\odel{o_{\delta}}
\def\od1{o_{\delta}(1)}

\def\oe1{o_{\varepsilon}(1)}
\def\aF{\alpha_F}
\def\bF{\beta_F}
\def\gF{\gamma_F}


\title{\bf Counting substructures I:  color critical graphs}
\author{Dhruv Mubayi
\thanks{ Department of Mathematics, Statistics, and Computer
Science, University of Illinois, Chicago, IL 60607.  email:
mubayi@math.uic.edu; research  supported in part by  NSF grant DMS
0653946
}}
\date{\today}
\maketitle

\begin{abstract}
Let $F$ be a  graph which contains an edge whose deletion reduces its chromatic number. We prove  tight bounds on the number of copies of $F$ in a graph with a prescribed number of vertices and edges. Our results extend those  of Simonovits \cite{S}, who proved that there is one copy of $F$, and of Rademacher, Erd\H os \cite{E1, E2} and Lov\'asz-Simonovits \cite{LS}, who proved similar counting results when $F$ is a complete graph.

One of the simplest cases of our theorem is the following new result. There is an absolute positive constant $c$ such that if $n$ is sufficiently large and $1 \le q <cn$, then every $n$ vertex graph with $\lfloor n^2/4 \rfloor+q$ edges contains at least
$$q\left \lfloor \frac{n}{2}\right\rfloor \left(\left \lfloor \frac{n}{2}\right\rfloor-1\right) \left(
\left\lceil \frac{n}{2}\right\rceil-2\right)$$
copies of a five cycle. Similar statements hold for any odd cycle and the bounds are best possible.

\end{abstract}

\section{Introduction}
 Mantel \cite{M} proved that  a graph with $n$ vertices and
$\lfloor n^2/4 \rfloor+1$ edges contains a triangle. Rademacher
extended this by showing that there are   at least $\lfloor n/2
\rfloor$ copies of a triangle.  Subsequently, Erd\H os \cite{E1, E2}
proved that if $q <cn$ for some small constant $c$, then $\lfloor
n^2/4 \rfloor+q$ edges guarantees at least $q\lfloor n/2 \rfloor$
triangles.  Later Lov\'asz and Simonovits \cite{LS} proved that the
same statement holds with $c=1/2$, thus confirming an old conjecture
of Erd\H os.  They also proved similar results for complete graphs.

In this paper we  extend the results of Erd\H os and Lov\'asz-Simonovits by
 proving such statements for the broader class of color critical graphs, which are graphs that contain an edge whose removal reduces their chromatic number\footnote{Note that a color critical graph is often defined as one with all proper subgraphs having lower chromatic number, but our definition is slightly  different}. In many ways our proof is independent of the specific structure of $F$.

The main new tool we use is the graph removal lemma, which is a
consequence of the hypergraph regularity lemma
 (see Gowers \cite{G}, Nagle-R\"odl-Schacht \cite{NRS}, R\"odl-Skokan \cite{RS}, Tao \cite{T}).  In subsequent papers \cite{M2, M3} we will extend these results to hypergraphs.
 The novelty in this project is the use of the removal lemma to count
 substructures in (hyper)graphs rather precisely.

We often associate a graph with its edge set. Given graphs
$F, H$, where  $F$ has $f$ vertices, a copy of $F$ in $H$ is a subset of
$f$ vertices and $|F|$ edges of $H$ such that the subgraph formed
by this set of vertices and edges is isomorphic to $F$.  In other words, if we denote $Aut(F)$ to be the number of automorphisms of $F$, then the number of copies of $F$ in $H$ is the number of edge-preserving injections from $V(F)$ to $V(H)$ divided by $Aut(F)$.
\medskip

\begin{theorem} {\bf (Graph Removal Lemma \cite{G, NRS, RS, T})} \label{removal}
Let $F$ be a graph with $f$ vertices. Suppose that an
$n$ vertex graph $H$ has at most $o(n^f)$ copies of $F$. Then
there is a set of edges in $H$ of size $o(n^2)$ whose removal from
$H$ results in a graph with no copies of $F$.
\end{theorem}

Say that a graph $F$ is $r$-critical if it has chromatic number $r+1$ and it contains an edge whose deletion reduces the chromatic number to $r$.  As usual, we define the Tur\'an number ex$(n,F)$ to be the maximum number of edges in an $n$ vertex graph that contains no copy of $F$ as a (not necessarily induced) subgraph.  The Tur\'an graph $T_r(n)$ is the $n$ vertex $r$-partite graph with the maximum number of edges; its parts all have size $\lceil n/r\rceil$ or $\lfloor n/r \rfloor$. Let
$$t_r(n)=|T_r(n)|=\sum_{1\le i<j\le r} \left\lfloor\frac{n+i-1}{r}\right\rfloor  \left\lfloor\frac{n+j-1}{r}\right\rfloor.$$  Since $\chi(T_r(n))=r$ it does not contain any $r$-critical graph $F$, and consequently ex$(n, F) \ge t_r(n)$.  Simonovits \cite{S} proved that if $n$ is sufficiently large, then we have equality.  In other words, every $n$ vertex graph ($n>n_0$) with $t_r(n)+1$ edges contains at least one copy of $F$. We extend his result by proving that there are  many copies and determining this optimal number. In fact, the number of copies is the number one gets by adding an edge to $T_r(n)$.

\begin{definition} Fix $r \ge 2$ and let $F$ be an $r$-critical graph. Then $c(n,F)$ is the minimum number of copies of $F$ in the graph obtained from $T_r(n)$ by adding one edge.
\end{definition}

Note that for any fixed $F$, computing $c(n, F)$ is just a finite process (see Lemma \ref{ab} in the next section). Our result below therefore gives an explicit formula for each color critical $F$, even though this formula may be very complicated.

\begin{theorem} \label{crit}
Fix $r \ge 2$ and an $r$-critical graph $F$. There exists  $\delta=\delta_F>0$ such that if $n$ is sufficiently large and  $1 \le q < \delta n$, then every $n$ vertex graph with $t_r(n)+q$ edges contains at least $q\, c(n,F)$ copies of $F$.
\end{theorem}

 Theorem \ref{crit} is asymptotically sharp, in that for every  $q<\delta n$, there exist graphs with $t_r(n)+q$ edges and at most $(1+O(1/n))q\, c(n,F)$ copies of $F$. To see this, simply add a matching of size $q$ to the appropriate part of $T_r(n)$. Each new edge lies in precisely $c(n,F)$ copies of $F$ that contain only one new edge, giving a total of $q\, c(n,F)$ copies.  If $F$ has $f$ vertices, then the number of copies of $F$ that contain at least two new edges is at most $O(q^2n^{f-4})$.  It is easy to see that $c(n,F)=\Theta(n^{f-2})$ (see Lemma \ref{ab} in the next section for more details), and since $q<n$, we obtain $q^2n^{f-4}=O(1/n)q\, c(n,F)$ as desired.

 In many instances Theorem \ref{crit} is sharp.  Let us examine two special cases.

{\bf Odd cycles.} Fix $k \ge 1$ and let $F=C_{2k+1}$. Then we quickly see that
$$c(n, F)=\lfloor n/2\rfloor(\lfloor n/2 \rfloor -1)\cdots(\lfloor n/2 \rfloor-k+1)
(\lceil n/2 \rceil-2)\cdots(\lceil n/2\rceil-k).$$
where we interpret the second product as empty if $k=1$. Moreover, if we add a matching within one of the parts of $T_2(n)$, then it is easy to see that no copy of $F$ contains two edges of the matching.  Hence Theorem \ref{crit} is sharp in the case $F=C_{2k+1}$. Even the simple case of counting $C_5$'s was not previously known.

{\bf $K_4$ minus an edge.} Let $F$ be the graph obtained from $K_4$ by deleting an edge. Then it is easy to see that $c(n, F)={\lfloor n/2 \rfloor \choose 2}$ and again this is sharp by adding a matching to one part.

For any fixed $\e>0$ and $1\le q < n^{1-\e}$, the proof of Theorem \ref{crit} actually produces a vertex that lies in $q \,c(n, F)$ copies of $F$ or $(1-\e)q$ edges that each lie in $(1-\e)c(n, F)$ copies of $F$. Such information about the distribution of the copies of $F$ does not seem to follow from the methods of \cite{E1, E2, LS}, even for cliques.

Throughout the paper, Roman alphabets (e.g. $r,s,n$) denote integers and Greek alphabets (e.g. $\aF, \bF, \gF, \e, \delta$) denote reals.  Given a set of pairs $H$, let $d_H(v)$ be the number of pairs in $H$ containing $v$.  So if we view $H$ as a graph, then $d_H(v)$ is just the degree of vertex $v$.

\section{Three lemmas}

In this section we will prove three technical lemmas needed in the proof of Theorem \ref{crit}.

\begin{lemma} \label{lemma1} Suppose that $r \ge 2$ is fixed, $n$ is sufficiently large, $s<n$ and $n_1+\cdots +n_r=n$. If
$$\sum_{1\le i<j \le r}n_in_j \ge t_r(n)-s,$$
then $\lfloor n/r \rfloor -s \le n_i \le \lceil n/r \rceil+s$ for all $i$.
\end{lemma}

\proof  The result is certainly true for $s=0$ by definition of $t_r(n)$ and the easy fact that $\sum_{1\le i<j \le r}n_in_j$ is maximized when $|n_i-n_j|<1$ for all $i$.  So assume that $s \ge 1$ and let us proceed by induction on $s$.  It is more convenient to prove the contrapositive, so assume that for some $i$, either $n_i>\lceil n/r\rceil+s$ or $n_i<\lfloor n/r \rfloor-s$ and we wish to prove that $\sum_{1\le i<j \le r}n_in_j < t_r(n)-s$.  We may assume that $n_1 \ge  \cdots \ge n_r$ and that
$n_1>\lceil n/r\rceil+s$ (the case  $n_r<\lfloor n/r \rfloor-s$ is symmetrical and has an almost identical proof).  Define $n_1'=n_1-1, n_r'=n_r+1$ and $n_i'=n_i$ for $1<i<r$. Then $\sum n_i'=n$ and we certainly have $n_i'>\lceil n/r\rceil+(s-1)$ for some $i$.  By the induction hypothesis, $$\sum_{1\le i<j \le r}n_i'n_j' < t_r(n)-(s-1).$$ We also have $n_1 \ge \lceil n/r\rceil+s+1$ and $n_r \le \lfloor n/r \rfloor$. Consequently,
\begin{align}
 \sum_{1\le i<j \le r}n_in_j&=\sum_{1\le i<j \le r}n_i'n_j'+n_r-(n_1-1)\notag \\
 &\le\sum_{1\le i<j \le r}n_i'n_j'+\lfloor n/r \rfloor-(\lceil n/r\rceil+s)\notag \\
 &\le
\sum_{1\le i<j \le r}n_i'n_j'-s\notag \\
&<t_r(n)-(s-1)-s\notag \\
&\le t_r(n)-s.\notag
\end{align}
This completes the proof of the lemma.
\qed
\medskip

\begin{lemma} \label{ab}
 Fix $r \ge 2$ and an $r$-critical graph $F$ with $f$ vertices. There are positive constants $\aF, \bF$ such that if $n$ is sufficiently large, then
 $$|c(n,F)-\aF n^{f-2}|<\bF n^{f-3}.$$
 In particular, $(\aF/2)n^{f-2}<c(n, F)< 2\aF n^{f-2}$.
 \end{lemma}
 \proof We may assume that $r|n$.  Indeed, suppose we could prove the result in the case $r|n$ and we are given $r$ that does not divide $n$.  The graph $T_r(n)$ has parts of size $\lfloor n/r \rfloor$ or $\lceil n/r \rceil$, so we can add at most one vertex to $r-1$ parts so that all parts have size $\lceil n/r\rceil$.  Also add edges from the new vertices to all parts distinct from the one that they lie in. The resulting  graph is $T_r(n')$ with  $n<n'<n+r$ vertices and $r|n'$. So we have
 $$c(n,F)\le c(n',F)\le c(n+r, F)<\aF(n+r)^{f-2}+\bF(n+r)^{f-3}<\aF n^{f-2}+(fr\aF+2\bF)n^{f-3}$$
 where the last inequality follows since $n$ is large.
 The theorem therefore holds with $\aF$ and $\beta'_F=fr\aF+2\bF$.  A similar argument gives the required lower bound on $c(n,F)$.

 Let us write an explicit formula for $c(n, F)$. Let $H$ be obtained from $T_r(n)$ by adding one edge $xy$ in the first  part. Say that an edge $uv \in F$ is good if $\chi(F-uv)=r$. Let $\chi_{uv}$ be a proper $r$-coloring of  $F-uv$  such that $\chi_{uv}(u)=\chi_{uv}(v)=1$.  Every proper $r$-coloring of $F-uv$ gives the same color to $u,v$, since $\chi(F)>r$. Let $x_{uv}^i$ be the number of vertices of $F$ excluding $u,v$ that receive color $i$. An edge preserving injection of $F$ to $H$ is obtained by choosing a good edge $uv$ of $H$, mapping it to $xy$, then mapping the remaining vertices of $F$ to $H$ such that no two adjacent vertices get mapped to the same part of $H$.  Such a mapping is given by a coloring $\chi_{uv}$, and the number of mappings associated with $\chi_{uv}$ is just the number of ways the vertices colored $i$ can get mapped to the $i$th part of $H$. The vertices mapped to the first part cannot get mapped to $x,y$ since these have been taken by $u,v$ and there are two ways to map $\{u,v\}$ to $\{x,y\}$. Altogether we obtain
 $$c(n,F)=\frac{1}{2^{f^2}}\sum_{uv\, good} \sum_{\chi_{uv}}2(n/r-2)_{x^1_{uv}}\prod_{i=2}^r(n/r)_{x^i_{uv}}.$$
 Expanding this expression, we get a sum of  polynomials in $n$ of degree $x^1_{uv}+\sum_{i=2}^r x^i_{uv}=f-2$ with positive leading coefficients. Consequently, $c(n,F)$ is a polynomial in $n$ of degree $f-2$ with positive leading coefficient.  Let $\aF$ be this coefficient.  Since $n$ is sufficiently large, and the coefficients of $c(n, F)$ are all fixed independent of $n$, we can bound the absolute value of the contribution of all the other terms by $\bF n^{f-3}$ for some positive $\bF$ that depends only on $F$.  This completes the proof. \qed
 \medskip

 Given integers $n_i, \ldots, n_r$, let $c(n_1, \ldots, n_r, F)$ be the number of copies of $F$   in the graph obtained from the complete $r$-partite graph with parts $n_1, \ldots, n_r$ by adding an edge to the part of size $n_1$.

\smallskip

\begin{lemma} \label{lemma2} Let $F$ be an $r$-critical graph. There is a positive constant $\gF$ depending only on $F$ such that the following holds for  $n$  sufficiently large. If $n_1+\cdots+n_r=n$ with $\lfloor n/r \rfloor -s\le n_i\le \lceil n/r\rceil + s$ and $s<n/3r$, then
$$c(n_1, \ldots, n_r, F)\ge c(n,F)-\gF s n^{f-3}.$$
\end{lemma}
\proof  If $s=0$, then the result holds by the definition of $c(n,F)$, so assume that $s \ge 1$. Let $H$ be the graph obtained from the complete $r$-partite graph with parts $n_1\ge  \ldots\ge n_r$ by adding an edge $xy$ to the part of size $n_1$. Since $|n_1-n_r|\le 2s+1$, we can remove at most $2s+1\le 3s$ vertices (excluding $xy$) from each part of $H$ so that each part has  size $n_r$.  The resulting graph $H'$ satisfies $H'-xy \cong T_r(n')$ with $n'\ge n-3rs$.  Since $n$ is sufficiently large and $s<n/3r$, we have
\begin{align}
(n-3rs)^{f-2}&=n^{f-2}+\sum_{i=0}^{f-3}{f-2 \choose i}n^i(-3rs)^{f-2-i} \notag \\
&=n^{f-2}-\sum_{i=0}^{f-3}{f-2 \choose i}n^i(-1)^{f-1-i}(3rs)^{f-2-i} \notag \\
&\ge n^{f-2}-\sum_{i=0}^{f-3}{f-2 \choose i}n^i(3rs)^{f-2-i} \notag \\
&=n^{f-2}-\sum_{i=0}^{f-3}{f-2 \choose i}(3r)^{f-2-i}sn^i s^{f-3-i} \notag \\
&>n^{f-2}-\sum_{i=0}^{f-3}{f-2 \choose i}\frac{(3r)^{f-2-i}}{(3r)^{f-3-i}}sn^{f-3} \notag \\
&>n^{f-2}-3rsn^{f-3}\sum_{i=0}^{f-3}{f-2 \choose i}\notag \\
&>n^{f-2}-2^fr s n^{f-3}. \notag
\end{align} Put $\gF=\aF 2^fr+2\bF$.
Lemma \ref{ab} now gives
\begin{align}
c(n_1, \ldots, n_r, F)&\ge c(n-3rs,F)
\notag \\
&>\aF(n-3rs)^{f-2}-\bF(n-3rs)^{f-3} \notag \\
&>\aF n^{f-2}-\aF 2^frsn^{f-3}-\bF n^{f-3} \notag \\
&\ge\aF n^{f-2} +\bF n^{f-3} -(\aF 2^fr+2\bF)sn^{f-3}\notag \\
&=\aF n^{f-2}+\bF n^{f-3}-\gF sn^{f-3} \notag \\
&>c(n,F)-\gF s n^{f-3} \notag
\end{align}
and the proof is complete.
\qed

\section{Proof of Theorem \ref{crit}}

In this section we will prove Theorem \ref{crit}.
We need the following stability result proved by Erd\H os and Simonovits \cite{S}.
\medskip

\begin{theorem}  {\bf (Erd\H os-Simonovits Stability Theorem \cite{S})} \label{stability}
Let $r \ge 2$ and $F$ be a fixed $r$-critical graph. Let $H$ be a graph with $n$ vertices and $t_r(n)-o(n^2)$ edges
that contains no copy of $F$. Then there is a partition of the
vertex set of $H$ into $r$ parts so that the number of edges contained within a part is at most $o(n^2)$. In other words, $H$ can
be obtained from $T_r(n)$ by adding and deleting a set of $o(n^2)$
edges.
\end{theorem}

{\bf Remark.} The $o(1)$ notation above should be interpreted in the obvious way, namely $\forall \eta,  \exists \xi, n_0$ such that if $n>n_0$ and $|H| >t_r(n)-\xi n^2$, then $H=T_r(n) \pm \eta n^2$ edges.  We will
not explicitly mention the role of $\xi, \eta$ when we use the result, but it should be obvious from the context. We will also assume that $\xi$ is much smaller than $\eta$, since if the result holds for $\xi$, then it also holds for $\xi'<\xi$.
Similar comments apply for  applications of Theorem \ref{removal}.

 \bigskip

\noindent{\bf Proof of Theorem \ref{crit}.}
Our constants $\delta_i, \e_i, \e$ will enjoy the following hierarchy:
$$1/n \lll \delta_1 \lll \delta_2 \ll \delta_3 \ll \delta_4 \ll \e_4 \ll \e_3 \ll \e_2 \ll \e_1 \ll \e \ll 1.$$ The notation $\xi \lll \eta$  above means that we are applying some theorem with input $\eta$ and output $\xi$ (the quantification is $\forall \eta, \exists \xi$)  and $\xi \ll \eta$ simply means that $\xi$ is a sufficiently small function of $\eta$ that is needed to satisfy some inequality in the proof.
In addition, we require
 $$\e<\frac{\aF}{4(\gF + 2^{f^2})}\,,$$ where $\aF$ comes from Lemma \ref{ab} and
 $\gF$ comes from Lemma \ref{lemma2}. Also, $n$ is sufficiently large that Lemmas \ref{lemma1} ,\ref{ab},  \ref{lemma2} apply whenever needed. We emphasize that $\e_4$ is an absolute constant that depends only on $F$. Set $\delta=\e_4/4\aF$ and suppose that $1 \le q<\delta n$. Let
 $H$ be an $n$ vertex graph
with $t_r(n)+q$ edges.  Write $\# F$ for the number of copies of
${F}$ in $H$.

If $\#{F}\ge n^{f-1/2}$, then since $c(n, F)<2\aF n^{f-2}$ and $q<(\e_4/4\aF) n$, we have $$\#F \ge n^{f-1/2}>(\e_4/4\aF)n (2\aF n^{f-2})> q\,c(n, F)$$
 and we are done. So assume that
$\#{F}<n^{f-1/2}=(1/n^{1/2}) n^f$.  Since $n$ is sufficiently large, by the Removal lemma there is a set of at
most $\delta_1 n^2$ edges of $H$ whose removal results in a
graph ${H}'$ with no copies of ${F}$. Since
$|{H}|>t_r(n)-\delta_1 n^2$, by Theorem \ref{stability},
we conclude that there is an $r$-partition of ${H}'$ (and also of
$H$) such that the number of edges contained entirely within a
part is $\delta_2n^2$. Now pick an $r$-partition $V_1 \cup \ldots \cup V_r$ of
$H$ that maximizes $e(V_1, \ldots, V_r)$, the number of edges that intersect
two  parts. We know that $e(V_1, \ldots, V_r)\ge t_r(n)-\delta_2n^2$, and an
easy calculation also shows that each  $V_i$ has size
$n_i=n/r\pm \delta_3n$.

Let $B$ (bad) be the set of edges of $H$ that lie entirely within a part and let $G$ (good) be the set of edges of $H$ that intersect two parts, so $G=H-B$. Let $M$ (missing) be the set of
pairs which intersect two parts that are not edges of $H$. Then
$G \cup M$ is  $r$-partite so it has at most $t_r(n)$ pairs and
$$t_r(n)+q-|B|+|M|=|(H-B) \cup M|=|G \cup M|=|G|+|M|\le t_r(n).$$
 Consequently,
$$q+|M| \le |B| \le \delta_2n^2.$$
Also, $|H|=|G|+|B|$ so we may suppose that $|G|=t_r(n)-s$ and $|B|=q+s$ for some $s \ge 0$.
For an edge $e \in B$, let $F(e)$ be the number of copies of $F$ in $H$ containing the unique edge $e$ from $B$.

If $s=0$, then $G \cong T_r(n)$ and $F(e)\ge c(n, F)$ for every $e \in B$ (by definition of $c(n, F)$) so we immediately obtain $\#F \ge |B|c(n,F)=qc(n, F)$.

 We may therefore assume that $s\ge 1$.  Partition $B=B_1 \cup B_2$, where
$$B_1=\{e \in B: F(e)>(1-\e)c(n,F)\}.$$
A potential copy of $F$ is a copy of $F$ in $G \cup M \cup B=H \cup M$ that uses exactly one edge of $B$.

\medskip

{\bf Claim 1.} $|B_1| \ge (1-\e)|B|$

Proof of Claim. Suppose to the contrary that $|B_2|\ge \e|B|$.
Pick $e \in B_2$. Assume wlog that $e \subset V_1$. By Lemma \ref{lemma2} (observing that $n_i =n/r \pm \delta_3n$ and $\delta_3n<n/3r$), the number of potential copies of $F$ containing $e$ is
$$c(n_1, \ldots, n_r, F) \ge c(n, F)-\gF (\delta_3 n) n^{f-3}>c(n,F)-\gF \delta_3n^{f-2}>(1-\delta_4)c(n, F).$$
 At
least $(\e/2)c(n,F)$ of these potential copies of ${F}$  have a
pair from $M$, for otherwise
$$F(e)>c(n_1, \ldots, n_r, F)-(\e/2)c(n, F)> (1-\delta_4-\e/2)c(n,{F})>(1-\e)c(n,{F})$$ which contradicts the definition of $B_2$.   Suppose that for at least $(\e/4)c(n,F)$ of these potential copies of ${F}$, the
pair from $M$ misses $e$.
The number of times each such pair from $M$ is counted is at most the number of ways to choose the remaining $f-4$ vertices of the potential copy of $F$ ($e$ is fixed), and then $|F|<f^2$ pairs among the chosen vertices so that the resulting vertices and pairs form a graph isomorphic to $F$.  There are
 at most $2^{f^2}n^{f-4}$ ways to do this.  We obtain the contradiction
$$\frac{\e \aF}{82^{f^2}} n^2=\frac{(\e/8)\aF n^{f-2}}{2^{f^2} n^{f-4}}<
\frac{(\e/4)c(n,F)}{2^{f^2}n^{f-4}}\le |M|<\delta_2n^2,$$
where the first strict inequality follows from  Lemma \ref{ab}.
We may therefore assume that for at least $(\e/4)c(n,F)$ of these potential copies of ${F}$, the
pair from $M$ intersects $e$.  Each such pair is counted at most $2^{f^2}n^{f-3}$ times, so we conclude that there exists $x \in e$ with
$$d_M(x)\ge \frac{(\e/8)c(n,F)}{2^{f^2}n^{f-3}}>\frac{(\e/8)(\aF/2)n^{f-2}}{2^{f^2} n^{f-3}}=\frac{\e\ \aF}{162^{f^2}} n> \e_1n.$$

Let
$$A=\{v \in V(H): d_M(v)>\e_1n\}.$$
We have argued above that every $e \in B_2$ has a vertex in $A$.
Consequently,
$$2\sum_{v \in A}d_{B_2}(v) \ge 2|B_2|\ge 2\e|B|> 2\e|M|\ge \e\sum_{v
\in A} d_M(v) >\e|A|\e_1 n,$$ and there exists a vertex $u \in  A$  such that
$d_{B_2}(u) \ge (\e\e_1/2)n>\e_2n$. Assume wlog that $u \in V_1$.

 By the choice of the partition, we may assume that $u$ has at least $\e_2n$ neighbors in  $V_i$ for each $i>1$, otherwise moving $u$ to $V_i$ increases $e(V_1, \ldots, V_r)$.  For each $i=1, \ldots, r$, let $V_i'$ be a set of $\lceil \e_2n\rceil $ neighbors of $u$ in $V_i$ and for $i=1$, enlarge $V_1'$ by  letting $u \in V_1'$.  Let $n_i'=|V_i'|$.  Pick $v \in V_1'-\{u\}$. Then the number of  copies of $F$ in the complete $r$-partite graph $K(V_1', \ldots, V_r')$ together with edge $e=uv$ is by definition at least
 $$c(r\lceil \e_2n\rceil, F)>\aF(r\e_2n)^{f-2}-\bF(r\lceil \e_2n\rceil)^{f-3}>\e_3n^{f-2}$$
 for suitable $\e_3$ depending only on $F$. Let us sum this inequality over all such $e=uv \in B_2$ with $v \in V_1'$. Since $\delta=\e_4/4\aF$ and $c(n,F)<2\aF n^{f-2}$, we obtain at least $$(|V_1'|-1)\e_3 n^{f-2}\ge (\e_2n)\e_3n^{f-2}>\e_4n^{f-1}=4\delta \aF n^{f-1}>4q\aF n^{f-2}>2q\, c(n,F)$$ potential copies of $F$ containing $u$.
 At least half of these potential copies of $F$ must have a pair from $M$, otherwise we are done.  This pair from $M$ cannot be incident with $u$, since $u$ is adjacent to all vertices in $K(V_1', \ldots, V_r')$ (other than itself).  Hence this pair from $M$ misses $u$. Each such pair is counted at most $2^{f^2}n^{f-3}$ times, so we obtain the contradiction
$$\frac{\e_4}{ 2^{f^2+1}} n^2=\frac{\e_4n^{f-1}}{2^{f^2+1}n^{f-3}}<|M|<\delta_2n^2.$$ This concludes the proof of the Claim. \qed
\bigskip

If $s\ge 4\e q$, then counting copies of $F$ from edges of $B_1$ and using Claim 1 we get
\begin{align}
\#F \ge \sum_{e \in B_1} F(e)&\ge \sum_{e \in B_1}(1-\e)c(n,F)\notag \\
&\ge |B_1|(1-\e)c(n,F) \notag \\
&\ge (1-\e)^2|B|c(n, F)\notag \\
&>(1-2\e)(q+s)c(n,F)\notag \\
&\ge (q+2\e q-8\e^2q)c(n,F)\notag \\
&> q\, c(n,F).\notag
\end{align}
So we may assume that $s <4\e q<q<n/3r$. Recall that $n_i=|V_i|$.
Then
$$t_r(n)-s=|G| \le \sum_{1\le i<j \le r}n_in_j.$$
Now  Lemma \ref{lemma1} implies that for each $i$,
$$\lfloor n/r \rfloor -s\le n_i\le \lceil n/r\rceil + s.$$

Observe that $|M| \le s$ for otherwise $|G\cup M|>t_r(n)$ which is impossible.
Pick $e \in B$ and assume wlog that $e \subset V_1$.
The number of potential copies of $F$ containing $e$ is by definition
$c(n_1, \ldots, n_r, F)$.
Now
 Lemma \ref{lemma2} implies that
 $$c(n_1, \ldots, n_r, F) \ge  c(n, F)-\gF sn^{f-3}.$$ Not all of these potential copies of $F$ are in $H$, in fact, a pair from $M$ lies in at most $2^{f^2}n^{f-3}$ potential copies counted above (we may assume that the pair intersects $e$ otherwise it is counted at most $2^{f^2}n^{f-4}$ times). We conclude that
\begin{equation} \label{fe}
F(e) \ge c(n_1, \ldots, n_r, F)-2^{f^2}n^{f-3}|M|\ge  c(n, F)-\gF sn^{f-3}-2^{f^2}sn^{f-3}.\notag\end{equation}
Since $s<q$ this implies that
\begin{align} \# F \ge \sum_{e \in B}F(e) &\ge (q+s)(c(n, F)-\gF sn^{f-3}-2^{f^2}sn^{f-3}) \notag \\
&> q\,c(n, F)+s\,c(n, F)-2q(\gF sn^{f-3}+2^{f^2}sn^{f-3}).\notag
\end{align}
As $s<q<\delta n < \e n$ and $\e<\aF/(4(\gF + 2^{f^2}))$, we have the bound
$$2q(\gF+2^{f^2})sn^{f-3} <2\e(\gF+2^{f^2})sn^{f-2}<s(\aF/2)n^{f-2}<s \,c(n,F)$$
This shows that $\#F>q\,c(n, F)$ and completes the proof of the theorem. \qed

\medskip
\section{Acknowledgments}
 I wish to thank Hemanshu Kaul who was instrumental in starting this project by posing a very general question about the number of  copies of certain subgraphs that one is guaranteed to find on a vertex. Initial discussions between Kaul and the author led to a proof of a weaker version of Theorem \ref{crit} that found $(1-o(1))q\,c(n, F)$ copies of $F$.

\end{document}